\documentclass[12pt]{article}
\usepackage{amssymb}
\usepackage{array}
\usepackage{latexsym}
\usepackage{enumerate}
\usepackage{amsmath}
\usepackage{amsfonts}
\usepackage{amsthm}
\usepackage{psfrag}
\usepackage{comment}
\usepackage[english]{babel}

\title{\bf{On large automorphism groups of algebraic curves
in positive characteristic}}
\author{M.~Giulietti ${}^*$ and G.~Korchm\'aros
 \thanks{Research supported by  the Italian
    Ministry MURST, Strutture geometriche, combinatoria e loro
    applicazioni, PRIN 2006-2007}}

\newtheorem{theorem}{Theorem}[section]

\newtheorem{lemma}[theorem]{Lemma}

{\theoremstyle{definition}
\newtheorem*{definition*}{Definition}

\newtheorem*{proposition*}{Proposition}
\newtheorem*{corollary*}{Corollary}
\newtheorem*{lemma*}{Lemma}

\def\cE{\mathcal E}
\def\cF{\mathcal F}

\def\cL{\mathcal L}

\def\cX{\mathcal X}
\def\cY{\mathcal Y}
\def\cZ{\mathcal Z}

\def\K{\mathbb{K}}
\def\PG{{\rm{PG}}}
\def\ord{\mbox{\rm ord}}
\def\deg{\mbox{\rm deg}}

\newcommand{\PSL}{\mbox{\rm PSL}}
\newcommand{\PGL}{\mbox{\rm PGL}}

\newcommand{\PSU}{\mbox{\rm PSU}}
\newcommand{\PGU}{\mbox{\rm PGU}}

\newcommand{\Sz}{\mbox{\rm Sz}}
\newcommand{\Ree}{\mbox{\rm Ree}}

\newcommand{\aut}{\mbox{\rm Aut}}

\newcommand{\ga}{\alpha}

\newcommand{\gS}{\Sigma}
\newcommand{\gO}{\Omega}

\newcommand{\ha}{{\textstyle\frac{1}{2}}}

\newcommand{\bv}{{\bf v}}

\sloppy

\begin{document}
\maketitle

    \begin{abstract}
In his investigation on large $K$-automorphism groups of an
algebraic curve, Stichtenoth obtained an upper bound on the order
of the first ramification group of an algebraic curve $\cX$
defined over an algebraically closed field of characteristic $p$.
Stichtenoth's bound has raised the problem of classifying all
$\K$-automorphism groups $G$ of $\cX$ with the following property:
There is a point $P\in \cX$ for which
\begin{equation}
\label{eq0} |G_P^{(1)}|>\frac{p}{p-1}g.
\end{equation}
Such a classification is obtained here by proving Theorem
\ref{refinementstichtenotharchivstaz1}
    \end{abstract}
    \section{Introduction} Let $\K$ be an algebraically
closed field of characteristic $p\geq 0$. Let $\aut(\cX)$ be the
$\K$-automorphism group of a (projective, non-singular,
geometrically irreducible, algebraic) curve $\cX$ embedded in an
$r$-dimensional projective space $\PG(r,\K)$. If $\cX$ has genus
$g\geq 2$, then $\aut(\cX)$ is finite, but this does not hold for
either rational or elliptic function fields. Over $\mathbb{C}$,
more generally in zero characteristic, the classical Hurwitz bound
for curves of genus $g\geq 2$ is
\begin{equation}
\label{hbound} |\aut(\cX)|\leq 84(g-1).
\end{equation}
In positive characteristic, $\cX$ with $g\geq 2$ may happen to
have larger automorphism group. Nevertheless, from previous work
by Roquette \cite{roquette1970}, Stichtenoth
\cite{stichtenoth1973I,stichtenoth1973II}, Henn \cite{henn1978},
Hansen and Petersen \cite{hansen1993}, Garcia, Stichtenoth and
Xing \cite{garcia-stichtenoth-xing2001}, \c{C}ak\c{c}ak and
\"Ozbudak,\cite{cakcak-ozbudak2004},  Giulietti, Korchm\'aros and
Torres \cite{giulietti-korchmaros-torres2004}, Lehr and Matignon
\cite{lehr-matignon2005}, curves with very large automorphism
groups comparing with their genera are somewhat rare.

The Hermitian curve is the unique curve with $|\aut(\cX)|\geq
16g^4$, see \cite{stichtenoth1973I}. If $p>g$ then (\ref{hbound})
holds with only one exception, namely the hyperelliptic curve
$\bv(Y^p-Y-X^2)$
with $g=\ha (p+1)$ and $|G|=2p(p^2-1)$,
see \cite{roquette1970}. Curves with $|\aut(\cX)|\geq 8g^3$ were
classified in \cite{henn1978}:
\begin{theorem}[Henn]
\label{hennmainbound} Let $\cX$ be a projective, geometrically
irreducible non-singular curve of genus $g\geq 2$. If a subgroup
$G$ of $|\aut(\cX)|$ has order at least $8g^3,$ then $\cX$ is
birationally equivalent to one of the following plane curves$:$
\begin{enumerate}
    \item[\rm(I)] The hyperelliptic curve $\bv(Y^2+Y +X^{2^k+1})$ with
$p=2$ and $g=2^{k-1};$ $|\aut(\cX)|=2^{2k+1} (2^k+1)$ and
$\aut(\cX)$ fixes a point $P\in\cX$.

    \item[\rm(II)] The hyperelliptic curve  $\bv(Y^2
-(X^q-X))$ with $p>2$ and $g=\ha(n-1);$  either $G/M\cong
\PSL(2,n)$\, or\, $G/M\cong \PGL(2,n),$ where $n$ is a power of
$p$ and $M$ is a central subgroup of $G$ of order $2$.
    \item[\rm(III)] The Hermitian curve  $\bv(Y^n+Y - X^{n+1})$ with
$p\geq 2$, and $g=\ha(n^2-n);$ either $G\cong \PSU(3,n)$\, or\,
$G\cong \PGU(3,n)$ with $n$ a power of $p$.

    \item[\rm(IV)] The {\rm{DLS}} curve (the Delign-Lusztig curve arising from the Suzuki group) $\bv(X^{n_0}(X^n+X) -
(Y^n+Y))$ with $p=2,\,n_0=2^r, r\geq 1, n=2n_0^2,$ and
$g=n_0(n-1);$ $G\cong \Sz(n)$ where $\Sz(n)$ is the Suzuki group.
\end{enumerate}
\end{theorem}
Another relevant example in this direction is the following.
\begin{enumerate}{\em
    \item[\rm(V)] The {\rm{DLR}} curve (the Delign-Lusztig curve arising from the Ree group) $\bv(X^{n_0}(X^n+X) -
(Y^n+Y))$ with $p=3,\,n_0=3^r,\ n\geq 0,\,n=3n_0^2,$ and
$g=n_0(n-1);$ $G\cong \Ree(n)$ where $\Ree(n)$ is the Ree group.}
\end{enumerate}
An important ingredient in the proof of Theorem
\ref{hennmainbound}, as well as in other investigations on curves
with a large automorphism group $G$, is Stichtenoth's upper bound
on the order of the Sylow $p$-subgroup in the stabiliser $G_P$ of
a point $P\in \cX$. In terms of ramification groups, such a Sylow
$p$-subgroup is the first ramification group $G_P^{(1)}$, and the
bound depends on the ramification pattern of the Galois covering
$\cX\to \cY$ where $\cY$ is a non-singular model of the quotient
curve $\cX/G_P^{(1)}$.
\begin{theorem}[Stichtenoth]
\label{stichtenotharchivsatz1} Let $\cX$ be a projective,
geometrically irreducible, non-singular curve of genus $g\geq 2$.
If $P\in\cX$, then
$$
|G_{P}^{(1)}|\leq \frac{4p}{p-1} g^2.
$$
More precisely$,$ if $\cX_i$ is the quotient curve
$\cX/G_{P}^{(i)}$, then one of the following cases occurs$:$
\begin{enumerate}

\item[\rm(i)] $\cX_1$ is not rational$,$ and $|G_{P}^{(1)}|\leq
g;$

\item[\rm(ii)] $\cX_1$ is rational, the covering $\cX\to\cX_1$
ramifies not only at $P$ but at some other point of $\cX$, and
\begin{equation}
\label{henning1}
 |G_P^{(1)}|\leq \frac{p}{p-1}g;
\end{equation}

\item[\rm(iii)] $\cX_1$ and  $\cX_2$ are rational$,$ the covering
$\cX\to\cX_1$ ramifies only at $P$, and
\begin{equation}
\label{henning2} |G_{P}^{(1)}| \leq
\frac{4|G_{P}^{(2)}|}{(|G_{P}^{(2)}|-1)^2}g^2\leq
\frac{4p}{(p-1)^2}g^2.
\end{equation}
\end{enumerate}
\end{theorem}
Stichtenoth's bound raises the problem of classifying all
automorphism groups $G$ with the following property: There is a
point $P\in \cX$ such that
\begin{equation}
\label{eq0bis} |G_P^{(1)}|>\frac{p}{p-1}g.
\end{equation}
In this paper, we obtain such a classification by proving the
following result.
\begin{theorem}
\label{refinementstichtenotharchivstaz1} If {\rm{(\ref{eq0})}}
holds, then either $G$ fixes $P$ or one of the four cases
{\rm{(II),$\ldots$,(V)}} in Theorem \ref{hennmainbound} occurs.
\end{theorem}

    \section{Background and some preliminary results}
    \label{background}
Let $\cX$ be a projective, non-singular, geometrically
irreducible, algebraic curve of genus $g\geq 2$ embedded in the
$r$-dimensional projective space $\PG(r,\K)$ over an algebraically
closed field $\K$ of positive characteristic $p>0$. Let $\gS$ be
the function field of $\cX$ which is an algebraic function field
of transcendency degree one over $\K$. The automorphism group
$\aut(\cX)$ of $\cX$ is defined to be the automorphism group of
$\gS$ fixing every element of $\K$. It has a faithful permutation
representation on the set of all points $\cX$ (equivalently on the
set of all places of $\gS$). The orbit
$$o(P)=\{Q\mid Q=P^\alpha,\, \alpha\in G\}$$
is {\em long} if $|o(P)|=|G|$, otherwise $o(P)$ is {\em short} and
$G_P$ is non-trivial.

If $G$ is a finite subgroup of $\aut(\cX)$, the subfield $\gS^G$
consisting of all elements of $\gS$ fixed by every element in $G$,
also has transcendency degree one. Let $\cY$ be a non-singular
model of $\gS^G$, that is, a projective, non-singular,
geometrically irreducible, algebraic curve with function field
$\gS^G$. Sometimes, $\cY$ is called the quotient curve of $\cX$ by
$G$ and denoted by $\cX/G$. The covering $\cX\mapsto \cY$ has
degree $|G|$ and the field extension $\gS/\gS^G$ is of Galois
type.

If $P$ is a point of $\cX$, the stabiliser $G_P$ of $P$ in $G$ is
the subgroup of $G$ consisting of all elements fixing $P$. For
$i=0,1,\ldots$, the $i$-th ramification group $G_P^{(i)}$ of $\cX$
at $P$ is
$$G_P^{(i)}=\{\alpha\mid \ord_P(\alpha(t)-t)\geq i+1, \alpha\in
G_P\}, $$ where $t$ is a uniformizing element (local parameter) at
$P$. Here $G_P^{(0)}=G_P$ and $G_P^{(1)}$ is the unique Sylow
$p$-subgroup of $G_P$.
Therefore, $G_P^{(1)}$ has a cyclic complement $H$ in $G_P$, that
is,
$G_P=G_P^{(1)}\rtimes H$ with a cyclic group $H$ of order prime to
$p$. Furthermore, for $i\geq 1$, $G_P^{(i)}$ is a normal subgroup
of $G$ and the factor group $G_P^{(i)}/G_P^{(i+1)}$ is an
elementary abelian $p$-group. For $i$ big enough, $G_P^{(i)}$ is
trivial.

For any point $Q$ of $\cX$, let $e_Q=|G_Q|$ and
$$d_Q=\sum_{i\geq 0}(|G_Q^{(i)}|-1).$$ Then $d_Q\geq e_Q-1$
and equality holds if and only if $\gcd(p,|G_Q|)=1.$

Let $g'$ be the genus of the quotient curve $\cX/G$. The Hurwitz
genus formula together with the Hilbert different formula give the
following equation
    \begin{equation}
    \label{eq1}
2g-2=|G|(2g'-2)+\sum_{Q\in \cX} d_Q.
    \end{equation}

If $G$ is tame, that is $p\nmid |G|$, or more generally for $G$
with $p\nmid e_Q$ for every $Q\in\cX$, Equation (\ref{eq1}) is
simpler and may be written as
    \begin{equation}
    \label{eq2}
2g-2=|G|(2g'-2)+\sum_{i=1}^k (|G|-|\ell_i|)
    \end{equation}
where $\ell_1,\ldots,\ell_k$ are the short orbits of $G$ on $\cX$.

Let $G_P=G_P^{(1)}\rtimes H$. The following upper bound on $|H|$
depending on $g$ is due to Stichtenoth \cite{stichtenoth1973I}:
\begin{equation*}
\label{stichtenotharchivsatz2} |H|\leq 4g+2.
\end{equation*}

For any abelian subgroup $G$ of $\aut(\cX)$, Nakajima
\cite{nakajima1987bis} proved that
\begin{equation*}
|G|\leq \left \{
\begin{array}{ll}
4g+4 & \mbox{for\quad  $p\neq 2,$} \\
4g+2 & \mbox{for\quad  $p=2.$}
\end{array}
\right.
\end{equation*}

Let $\cL$ be the projective line over $\K$. Then $\aut(\cL)\cong
\PGL(2,\K),$ and $\aut(\cL)$ acts on the set of all points of
$\cL$ as $\PGL(2,\K)$ naturally on $\PG(2,\K)$. In particular, the
identity of $\aut(\cL)$ is the only automorphism in $\aut(\cL)$
fixing at least three points of $\cL$. Every automorphism
$\ga\in\aut(\cL)$ fixes a point; more precisely, $\ga$ fixes
either one or two points according as its order is $p$ or
relatively prime to $p$. Also, $G_P^{(1)}$ is an infinite
elementary abelian $p$-group. For a classification of subgroups of
$\PGL(2,\K)$, see \cite{maddenevalentini1982}.

Let $\cE$ be an elliptic curve. Then $\aut(\cE)$ is infinite;
however for any point $P\in \cE$ the stabiliser of $P$ is rather
small, namely
\begin{equation*}
|\aut(\cE)_{P}| = \left \{
\begin{array}{ll}
 2,4,6 & \mbox{\quad when $p\neq 2,3,$} \\
 2,4,6,12 & \mbox{\quad when $p=3,$} \\
 2,4,6,8,12,24 & \mbox{\quad when $p=2.$}
\end{array}
\right.
\end{equation*}

Let $\cF$ be a (hyperelliptic) curve of genus $2$. For any
solvable subgroup $G$ of $\aut(\cF)$, Nakajima's bound together
with some elementary facts on finite permutation groups, yield
$|G|\leq 48$.

In the rest of this Section, $\cX$ stands for a non--hyperelliptic
curve of genus $g\geq 3$, and it is assumed to be the canonical
curve of $\K(\cX)$. So, $\cX$ is a non-singular curve of degree
$2g-2$ embedded in $\PG(g-1,\K)$, and the canonical series of
$\K(\cX)$ is cut out on $\cX$ by hyperplanes. Let
$1,x_1,\ldots,x_{g-1}$ denote the coordinate functions of this
embedding with respect to a homogeneous coordinate frame
$(X_0,X_1,\ldots,X_{g-1})$ in $\PG(g-1,\K)$.

For a point $P\in \cX$, the order sequence of $\cX$ at $P$ is the
strictly increasing sequence
\begin{equation}
\label{ordseq} j_0(P)=0<j_1(P)=1<j_2(P)<\ldots< j_{g-1}(P)
\end{equation}
such that each $j_i(P)$ is the intersection number $I(P,\cX\cap
H_i)$ of $\cX$ and some hyperplane $H_i$ at $P$, see \cite{sv}.
For $i=g-1$, such a hyperplane $H_{g-1}$ is uniquely determined
being the osculating hyperplane to $\cX$ at $P$. Another
characterisation of the integers $j_i(P)$, called $P$-orders or
Hermitian $P$-invariants, appearing in (\ref{ordseq}) is that $j$
is a $P$-order  if and only if $j+1$ is a Weierstrass gap, that
is, no element in $\K(\cX)$ regular outside $P$ has a pole of
order $j+1$. Now, assume that $j_{g-1}(P)=2g-2$, that is, $P$ is
the unique common point of $H_{g-1}$ with $\cX$. Then the
hyperplanes of $\PG(g-1,\K)$ whose intersection number with $\cX$
at $P$ is at least $j_{g-2}(P)$ cut out on $\cX$ a linear series
$g_n^1$ of degree $n=2g-2-j_{g-2}(P)$ and projective dimension
$1$. Let $\ell$ be the projective line over $\K$. Then $g_n^1$
gives rise to a covering $\cX\to\ell$ of degree $n$
which completely ramifies at $P$. If it also ramifies at the
points $P_1,\ldots,P_k$ of $\cX$ other than $P$, that is
$e_{P_i}>1$ for $i=1,\ldots,k$, then (\ref{eq1}) yields
\begin{equation}
\label{gnr} 2g-2=-2n+d_P+\sum_{i=1}^k d_{P_i}\geq
-(n+1)+\sum_{i=1}^k d_{P_i}.
\end{equation}
Note that $n$ must be at least $3$ as $\cX$ is neither rational,
nor elliptic and nor hyperelliptic.

{}From finite group theory, the following results and permutation
representations play a role in the proofs.

{\em Huppert's classification theorem}, see \cite[Chapter
XII]{huppert-blackburn1982III}: Let $G$ be a solvable
$2$-transitive permutation group of even degree $n$. Then $n$ is a
power of $2$, and $G$ is a subgroup of the affine semi-linear
group ${\mathrm{A}}\Gamma L(1,n)$.

{\em The Kantor-O'Nan-Seitz theorem}, see
\cite{kantor-o'nan-seitz1972}: Let $G$ be a finite $2$-transitive
permutation group whose $2$-point stabiliser is cyclic. Then $G$
has either a regular normal subgroup$,$ or $G$ is one of the
following groups in their natural $2$-transitive permutation
representations$:$
$$\PSL(2,n),\,\PGL(2,n),\,\PSU(3,n),\,\PGU(3,n),\Sz(n),\,\Ree(n).$$

{\em The natural $2$-transitive permutation representations of the
above linear groups:}
\begin{itemize}
\item[(i)] $G=\PGL(2,n)$, is the automorphism group of $\PG(1,n)$;
equivalently, $G$ acts on the set $\Omega$ of all
${\mathbb{F}}_{n}$-rational points of the projective line defined
over ${\mathbb{F}}_n$. The natural $2$-transitive representation
of $\PSL(2,n)$ is obtained when $\PSL(2,n)$ is viewed as a
subgroup of $\PGL(2,n)$,

\item[(ii)] $G=\PGU(3,n)$ is the linear collineation group
preserving the classical unital in the projective plane
$\PG(2,n^2)$, see \cite{hirschfeld1998}; equivalently $G$ is the
automorphism group of the Hermitian curve regarded as a plane
non-singular curve defined over the finite field ${\mathbb{F}}_n$
acting on the set $\Omega$ of all ${\mathbb{F}}_{n^2}$-rational
points. $\PSU(3,n)$ can be viewed as a subgroup of $\PGU(3,n)$ and
this is the natural $2$-transitive representation of $\PSU(3,n)$.

\item[(iii)] $G=\Sz(n)$ with $n=2n_0^2$, $n_0=2^r$ and $r\geq 1$,
is the linear collineation group of $\PG(3,n)$ preserving the Tits
ovoid, see \cite{tits1960,tits1962,hirschfeld1985}; equivalently
$G$ is the automorphism group of the {\rm{DLS}} curve regarded as
a non-singular curve defined over the finite field
${\mathbb{F}}_n$ acting on the set $\Omega$ of all
${\mathbb{F}}_{n}$-rational points.

\item[(iv)] $G=\Ree(n)$ with $n=3n_0^2$, $n_0=3^r$  and $r\geq 0$,
is the linear collineation group of $PG(7,n)$ preserving the Ree
ovoid, see \cite{tits1960}; equivalently, $G$ is the automorphism
group of the {\rm{DLS}} curve regarded as a non-singular curve
defined over the finite field ${\mathbb{F}}_n$ acting on the set
$\Omega$ of all ${\mathbb{F}}_{n}$-rational points.
\end{itemize}

For each of the above linear groups, the structure of the
$1$-point stabilizer and its action in the natural $2$-transitive
permutation representation, as well as its automorphism group, are
explicitly given in the papers quoted.

{\em Cyclic fix-point-free subgroups of some $2$-transitive
groups.} The following technical lemma is a corollary of the
classification of subgroups of $\PSU(3,n)$ and $\Ree(n)$.
\begin{lemma}
\label{backgr} Let $G$ be a $2$-transitive permutation group of
degree $n$. Let $U$ be a cyclic subgroup of $G$ which contains no
non-trivial element fixing a point.
    \begin{itemize}
    \item[(i)] If $G=\PSU(3,n)$ in its natural $2$-transitive
    permutation representation, then $|U|$ divides either $n+1$ or
    $n^2-n+1$.
    \item[(ii)] If $G=\Sz(n)$ in its natural $2$-transitive
    permutation representation, then $|U|$ divides either $n+1$,
    or $n-2n_0+1$, or $n+2n_0+1$.
    \item[(iii)] If $G=\Ree(n)$ in its natural $2$-transitive
    permutation representation, then $|U|$ divides either $n+1$,
    or $n-3n_0+1$, or $n+3n_0+1$.
    \end{itemize}
\end{lemma}
{\em Schur multiplier of some simple groups.} For a finite group
$G$, a group $\Gamma$ is said to be a covering of $G$ if $\Gamma$
has a central subgroup $U$, i.e. $U\subseteq Z(\Gamma)$, such that
$G\cong \Gamma/U$. If, in addition, $\Gamma$ is perfect, that is
$\Gamma$ coincides with its commutator subgroup, then  the
covering is called proper. For a simple group $G$, a perfect
covering is also called a semisimple group. From Schur's work, see
\cite{aschbacher1993} and \cite[V.23,24,25]{huppertI1967}, if $G$
is a simple group, then it possesses a ``universal'' proper
covering group $\bar{\Gamma}$ with the property that every proper
covering group of $G$ is a homomorphic image of $\bar{\Gamma}$.
The center $Z(\bar{\Gamma})$ is called the Schur multiplier of
$G$. The Schur multipliers of simple groups are known, see Griess
\cite{griess1972,griess1980,gorenstein1982}. In particular, the
Schur multiplier of $\PSL(2,q)$ with $q\geq 5$ odd, has order $2$;
$\PSU(3,q)$ with $q\geq 3$ has non-trivial Schur multiplier only
for $3|(q+1)$, and if this occurs the Schur multiplier has order
$3$; $\Ree(n)$ with $n>3$ has trivial Schur multiplier. Therefore,
the following result holds.
\begin{lemma}
\label{schur} Let $G$ be a simple group isomorphic to either
$\PSU(3,n)$ with $n\geq3$, or $\Ree(n)$ with $n>3$. If the center
$Z(\Gamma)$ of a group $\Gamma$ has order $2$ and $G\cong
\Gamma/Z(\Gamma)$ then $\Gamma$ has a subgroup isomorphic to $G$
and $\Gamma=Z(G)\times G$.
\end{lemma}

 \section{Large $p$-subgroups of $\aut(\cX)$ fixing a point}
In this section, Theorem \ref{refinementstichtenotharchivstaz1} is
proven.

We assume that case (iii) of Theorem \ref{stichtenotharchivsatz1}
with $G\neq G_P$ occurs.
In terms of the action of $G_P$ on $\cX$, Theorem
\ref{stichtenotharchivsatz1}
(iii)
implies that
\begin{itemize}
\item[(*)] {\em no non-trivial $p$-element in $G_P$ fixes a point
distinct from $P$.}
\end{itemize}
Let $\Omega$ be the set of all points $R\in\cX$ with non-trivial
first ramification group $G_R^{(1)}$. So, $\Omega$ consists of all
points $R\in \cX$ which are fixed by some element of $G$ of order
$p$. Since $P\in \Omega$ and $G\neq G_P$, $\Omega$ contains at
least two points. It may be noted that the $2$-point stabilizer of
$G$ is tame and hence cyclic.

Choose a non-trivial element $z$ from the centre of a Sylow
$p$-subgroup $S_p$ of $G$ containing $G_P^{(1)}$. Then $z$
commutes with a non-trivial element of $G_P^{(1)}$. This together
with (*) imply that $z$ fixes $P$. Therefore, $z\in G_P^{(1)}$. In
particular, $z$ fixes no point of $\cX$ distinct from $P$. Let
$g\in S_p$. Then $zg=gz$ implies that
$$(P^g)^z=P^{gz}=P^{zg}=(P^z)^g=P^g$$ whence $P^g=P$. This shows that
every element of $S_p$ must fix $P$, and hence $S_p=G_P^{(1)}$.
Since the Sylow $p$-subgroups of $G$ are conjugate under $G$,
every $p$-element fixes exactly one point of $\cX$. From Gleason's
Lemma, see \cite[Theorem 4.14]{dembowski1969}, $\Omega$ is a
$G$-orbit, and hence the unique non-tame $G$-orbit.

By (iii) of Theorem \ref{stichtenotharchivsatz1}, the quotient
curve $\cX_1=\cX/G_P^{(1)}$ is rational. This implies that
$\cY=\cX/G$ is also rational.

If there are at least two more short $G$-orbits, say $\Omega_1$
and $\Omega_2$, from (\ref{eq2}),
\begin{eqnarray*}
 2g-2 & \geq &-2|G|+(|G_P|+|G_P^{(1)}|+|G_P^{(2)}|-3)\,|\Omega|\\
&& \qquad + (|G_{Q_1}|-1)\,|\Omega_1|+(|G_{Q_2}|-1)\,|\Omega_2|,
\end{eqnarray*}
where $Q_i\in \Omega_i$ for $i=1,2$. Note that
$$
(|G_{Q_1}|-1)|\Omega_1|+(|G_{Q_2}|-1)|\Omega_2|\geq |G|
$$
since $|G_{Q_i}|-1\geq\ha |G_{Q_i}|$, and
$|G|=|G_{Q_i}|\,|\Omega_i|$. Also, $|G|=|G_{P}|\,|\Omega|$, and
$|G_{P}^{(2)}|>1$ by (iii) of Theorem
\ref{stichtenotharchivsatz1}. Therefore, $$ 2g-2\geq
(|G_{P}^{(1)}|-1)\, |\Omega|.
$$ Since
$|\Omega|\geq 2$, this implies that $g\geq |G_{P}^{(1)}|$, a
contradiction.



Therefore, one of the following cases occurs:
    \begin{itemize}
    \item[(i)] $\Omega$ is the unique short orbit of $G$;
    \item[(ii)] $G$ has two short orbits, namely $\Omega$ and a tame $G$-orbit.
    Furthermore, either
    \begin{itemize}
    \item[(iia)] there is a point $R\in \Omega$ such that the stabiliser of $R$ in $G_P$ is trivial;
    \item[(iib)] no point $R\in \Omega$ with the property as in (iia)
    exists.
    \end{itemize}
    \end{itemize}
Before investigating the above three cases separately, a useful
equation is established.

Since $\cX_1$ is rational, from (*)
\begin{equation}\label{formuladP1}
2g-2+2|G_P^{(1)}|=2|G_P^{(1)}|-2+|G_P^{(2)}|-1+\ldots=d_P-|G_P|+|G_P^{(1)}|.
\end{equation}
Therefore,
\begin{equation}\label{formuladP}
d_P=2g-2+|G_P^{(1)}|+|G_P|.
\end{equation}

\subsection{Case (i)}
\label{casei}
 We prove that only one example occurs, namely the DLR curve
for $n=3$.

{}From (*) it follows that $|\gO| \geq |G_P^{(1)}|+1$ whence
$|\gO|>g-1$ by (\ref{eq0}). On the other hand, (\ref{eq1}) gives
\begin{equation}
\label{eqrefinecaseiii} 2g-2=-2|G|+\deg D\,(\cX/G)= |\Omega|\,(
d_P
-2|G_{P}|).
\end{equation}
whence $|\Omega|$ is a divisor of $2g-2$. Therefore,
$|\Omega|=2g-2$ and $d_P=2|G_P|+1$. This and (\ref{formuladP})
give
\begin{equation}
\label{formuladP2}
 2g-2+|G_{P}^{(1)}|-|G_{P}|=1.
\end{equation}
  Since
$G_{P}=G_{P}^{(1)}\rtimes H$,
it follows that
$2g-2= (|H|-1)|G_{P}^{(1)}|+1$. This and (\ref{eq0}) imply that
$|H|=2$.

Therefore, $p\neq 2$ and $|G_{P}^{(1)}|=2g-3$. Since
$|\Omega|=2g-2$ and $G_{P}$ only fixes $P$, this implies that
$G_{P}^{(1)}$ acts on $\Omega\setminus \{P\}$ as a transitive
permutation group. Hence, $|\Omega|=q+1$ with $q=|G_{P}^{(1)}|$.
Since $\Omega$ is a $G$-orbit, it follows that $G$ induces on
$\Omega$ a $2$-transitive permutation group $\bar{G}$ whose
one-point stabiliser has order either $q$ or $2q$ according as
$|G|=2|\bar{G}|$ or $G=\bar{G}$.


If $|G|=2|\bar{G}|$, the subgroup $H$ is the kernel of the
permutation representation of $G$ on $\Omega$; that is, $H$ fixes
every point in $\Omega$. In particular, $H$ is a normal subgroup
of $G$. Therefore, $\bar{G}$ can be viewed as a $\K$--automorphism
group of the quotient curve $\cZ=\cX/G$. Let $P'$ be the point of
$\cZ$ lying under $P$. Since $2\nmid |G_{P}^{(1)}|$ it follows
that $G_{P}^{(1)}\cong \bar{G}_{P'}^{(1)}$. Also, the points of
$\cZ$ lying under the points in $\Omega$ form the unique short
$\bar{G}$-orbit. Therefore, case (i) occurs for $\cZ$ and
$\bar{G}$. {}From what was shown before, this implies that
$\bar{G}_{P'}=\bar{G}_{P'}^{(1)}\rtimes \bar{H}$ with
$|\bar{H}|=2$ But this is impossible as the stabiliser
$\bar{G}_{P'}$ has order $q$.

If $G=\bar{G}$, two cases are distinguished according as $G$ is
solvable or not. In the former case, Huppert's classification
theorem implies that $q+1=d^k$ with $d$ prime. Since $|\gO|=2g-2$
is even, so $d=2$. {}From Huppert's classification for $d=2$ it
also follows that $G_{P}$ is a subgroup of the $1$-point
stabiliser of $\mathrm{A\Gamma L}(1,q+1)$, and hence $|G_{P}|$
divides $kq$. On the other hand, $q+1=2^k$ can only occur when $k$
and $q$ are both primes. Since $|G_{P}|=2q$, this implies that
$k=2$. Hence $q=g=3$; that is, $p=3$ and $|G_{P}^{(1)}|=g=3$ which
contradicts (\ref{eq0}).

Suppose that $G$ is not solvable. If $G$ has a regular normal
subgroup $M$, then $G/M$ is not solvable. On the other hand,
$|G_P|=2q$ and $|M|=|\Omega|$. From $|G|=|G_P||\Omega|$, it
follows that $|G/M|=2q$. But this is not possible for a
non-solvable group, as $q$ is a prime power.

If $G$ does not have a regular normal subgroup, then we apply
Kantor-O'Nan-Seitz theorem. Since $|G|=2q$, this shows that either
$|\Omega|=6$ and $G\cong \PSL(2,5),$ or $|\Omega|=28$ and
$G\cong\Ree(3)$. In the former case, $|G_{P}^{(1)}|=5$ and $g=4$;
hence (\ref{eq0}) does not hold. In the latter case,
$|G_P^{(1)}|=27$ and $g=15$. This is consistent with (\ref{eq0}),
and $\cX$ is the smallest DLR curve. Therefore (V) holds for
$q=n^3$ with $n=3$.

\subsection{Case (iia)}
\label{caseiia}

We prove that no example occurs. Let $\Delta$ denote the unique
tame orbit of $G$. Choose a point $P$ from $ \Omega$ and a point
$Q$ from $\Delta$. Let
$$N=|G_Q|(d_P-|G_P|)-|G_P|.$$
Then
\begin{equation}
\label{Sticht13.7} |G|=2(g-1)\frac{|G_P^{(1)}||H||G_Q|}{N},
\end{equation}
where $G_P=G_P^{(1)}\rtimes H$.

By hypothesis, there exists a point $R\in \Omega$ such that the
the orbit $o(R)$ of $R$ under $G_P$ is long. Let $o'(R)$ denote
the orbit of $R$ under $G$. Then $|o'(R)|\cdot |G_{R}|=|G|$. Since
$P$ and $R$ lie in the same orbit $\Omega$ of $G$, so $G_{R}\cong
G_{P}$. Also, $o(R)$ is contained in $o'(R)$. Therefore,
\begin{equation}
\label{eq7} |G_{P}|\leq \frac{|G|}{|G_{P}|} = 2(g-1)\cdot
\frac{|G_{Q}|}{N}\leq 2(g-1)|G_{Q}|.
\end{equation}

Now, a lower bound on $N$ is given. As
$$
N\geq d_{P}|G_{Q}|- |G_{P}| |G_{Q}|-2(g-1)|G_{Q}|,
$$
so
\begin{equation}
\label{Sticht13.11} N\geq |G_{Q}|(d_{P}-|G_{P}|-2(g-1)).
\end{equation}
This and (\ref{formuladP}) imply that $N\geq
|G_{P}^{(1)}||G_{Q}|$. {}From (\ref{eq7}), $$N\leq
2(g-1)|G_Q|/|G_P|.$$ Hence $|G_P||G_P^{(1)}|\leq 2(g-1)$. Since
$|G_P|>1$, this contradicts (\ref{eq0}).

\subsection{Case (iib)}
\label{caseiib}
 We prove that $\cX$ is one of the examples
(II),(III),(IV) and (V) with $q>3$ in (V). Let $\Delta$ denote the
unique tame orbit of $G$. Choose a point $P$ from $ \Omega$ and a
point $Q$ from $\Delta$.

First the possible structure of $G$ and its action on $\Omega$ are
investigated.
\begin{lemma}
\label{2transitive} $G$ acts on $\Omega$ as a $2$-transitive
permutation group. In particular$,$ $|\Omega|=q+1$ with $q=p^t,$
and the possibilities for the permutation group $\bar{G}$ induced
by $G$ on $\Omega$ are as follows$:$
\begin{enumerate}
\item [\rm(1)] $\bar{G}\cong\PSL(2,q)$ or $\PGL(2,q);$

\item [\rm(2)] $\bar{G}\cong\PSU(3,n)$ or $\PGU(3,n),$ with
$q=n^3;$

\item [\rm(3)] $\bar{G}\cong \Sz(n),$ with $p=2,\ n=2n_0^2,\
n_0=2^k,$ with $k$ odd$,$ and $q=n^2;$

\item [\rm(4)]  $\bar{G}\cong \Ree(n)$ with $p=3,\ n=3n_0^2,\
n_0=3^k,$ and $q=n^3;$

\item [\rm(5)] a minimal normal subgroup of $\bar{G}$ is
solvable$,$ and the size of $\Omega$ is a prime power.
\end{enumerate}
\end{lemma}

\begin{proof} For a point $P\in \Omega$, let $\Omega_0=\{P\},\Omega_1,\ldots
\Omega_k$ with $k\geq 1$  denote the orbits of $G_{P}^{(1)}$
contained in $\Omega$. Then, $\Omega=\bigcup_{i=0}^k \Omega_i$. To
prove that $G$ acts $2$-transitively on $\Omega$,
it suffices to show that $k=1$.

For any $i$ with $1\leq i \leq k$, take a point $R\in \Omega_i$.
By hypothesis, $R$ is fixed by an element $\ga\in G_{P}$ whose
order $m$ is a prime different from $p$. Since
$|G_P|=|G_P^{(1)}||H|$ and $m$ divides $|G_P|$, this implies that
$m$ must divide $|H|$. By the Sylow theorem, there is a subgroup
$H'$ conjugate to $H$ in $G_{P}$ which contains $\ga$; here, $\ga$
preserves $\Omega_i$.

Since the quotient curve $\cX_1=\cX/G_P^{(1)}$ is rational, $\ga$
fixes at most two orbits of $G_{P}^{(1)}$. Therefore, $\Omega_0$
and $\Omega_i$ are the orbits preserved by $\ga$. As $H'$ is
abelian and $\ga\in H'$, this yields that $H'$ either preserves
both $\Omega_0$ and $\Omega_i$ or interchanges them. The latter
case cannot actually occur as $H'$ preserves $\Omega_0$. So, the
orbits $\Omega_0$ and $\Omega_i$ are also the only orbits of
$G_{P}^{(1)}$ which are fixed by $H'$. Since
$G_{P}=G_{P}^{(1)}\rtimes H'$, this implies that the whole group
$G_{P}$ fixes $\Omega_i$. As $i$ can be any integer between $1$
and $k$, it follows that $G_{P}$ fixes each of the orbits
$\Omega_0,\Omega_1,\ldots,\Omega_k$. Hence, either $k=1$ or
$G_{P}$ preserves at least three orbits of $G_{P}^{(1)}$. The
latter case cannot actually occur, as the quotient curve
$\cX_1=\cX/G_P^{(1)}$  is rational.

Therefore $k=1$. Also, the size of $\Omega$ is of the form $q+1$
with $q=|G_{P}^{(1)}|$; in particular, $q$ is a power of $p$.

Let $\bar{G}$ denote the $2$-transitive permutation group induced
by $G$ on $\Omega$. We apply the Kantor-O'Nan-Seitz theorem to
$\bar{G}$. Up to isomorphism, $\bar{G}$ is one of the groups on
the list, with $\bar{G}$ acting in each of the first four cases in
its natural $2$-transitive permutation representation.
\end{proof}
We also need the following consequence of Lemma \ref{2transitive}.
\begin{lemma}
\label{GQfpf} The subgroups $G_{P}$ and $G_{Q}$ have trivial
intersection, and $G_Q$ is a cyclic group whose order divides
$q+1$. Also,
\begin{equation}
\label{Sticht13.11tris}
2g-2=\frac{|G|\,(|G_{P}|-|G_{P}^{(1)}|\,|G_{Q}|)} {|G_{Q}|(|G| -
|G_{P}|)}
\end{equation}
\end{lemma}
\begin{proof} Let $\ga\in G_{P}\cap G_{Q}$ be non-trivial. Then
$p\nmid {\rm{ord}}\, \ga$, and hence $\ga\in H$. This shows that
$\ga$ fixes not only $P$ but another point in $\Omega$, say $R$.
Since $Q\not\in \Omega$, this shows that $\ga$ has at least three
fixed points. These points are in three different orbits of
$G_{P}^{(1)}$. Since the quotient curve $\cX_1=\cX/G_P^{(1)}$ is
rational, this implies that $\ga$ fixes every orbit of
$G_{P}^{(1)}$, a contradiction.
Hence $|G_{P}\cap G_{Q}|=1$. Therefore, no non-trivial element of
$G_Q$ fixes a point in $\Omega$. Since $|\Omega|=q+1$, the second
assertion follows.
Substituting $d_P$ from (\ref{formuladP}) into
$(\ref{Sticht13.7})$ gives (\ref{Sticht13.11tris}).
\end{proof}

First the case when the action of $G$ is faithful on $\Omega$ is
considered.

If $G\cong \PGL(2,q)$, then
$$
|G|=q^3-q,\ |G_{P}|=q^2-q,\ |G_{P}^{(1)}|=q.
$$
{}From Theorem \ref{stichtenotharchivsatz1}(iii), the second
ramification group $G_P^{(2)}$ is non-trivial.  As $G\cong
\PGL(2,q)$, $G_P$ has a unique conjugacy class of elements of
order $p$. Since $G_P^{(i)}$ is a normal subgroup of $G_P$, if
$g\in G_P^{(i)}$ with $i\geq 1$ then every conjugate of $g$ in $G$
also belongs to $G_P^{(i)}$. Therefore,
$$
G_P^{(1)}=G_P^{(2)}=\ldots= G_P^{(k)},\ |G_P^{(k+1)}|=1.$$ Since
the quotient curve $\cX_1=\cX/G_P^{(1)}$ is rational, from
(\ref{eq1}),
$$
2g=(q-1)(k-1).
$$ By (\ref{eq0}), this is only possible for $k=2$.
Therefore $g=\ha\,(q-1)$ with $q\geq 5$ odd, and
$|G_Q|=\ha\,(q+1)$.

Let $q\equiv 1 \pmod 4$. Then $2g-2\equiv 2 \pmod 4$, and
an involutory element in $\PGL(2,q)\setminus \PSL(2,q)$ has a no
fixed point on $\Omega$. Since $G_Q$ has odd order, such an
involutory element in $\PGL(2,q)\setminus \PSL(2,q)$ has no fixed
point in $\Delta$, either. Therefore, an involutory element in
$\PGL(2,q)\setminus \PSL(2,q)$ fixes no point of $\cX$. From
(\ref{eq2}) applied to such an involutory element
 $2g-2\equiv 0 \pmod 4$, a contradiction.

Let $q\equiv 3 \pmod 4$, then $2g-2\equiv 0 \pmod 4$
and an involutory element in $\PGL(2,q)\setminus \PSL(2,q)$ has
exactly two fixed points in $\Omega$. As before, $2\nmid |G_Q|$
implies that such an involutory element has no fixed point in
$\Delta$. Therefore, an involutory element in $\PGL(2,q)\setminus
\PSL(2,q)$ fixes exactly two points of $\cX$. From (\ref{eq2})
applied to such an involutory element
$2g-2\equiv 2\pmod 4$, a contradiction. Therefore, the case
$G\cong \PGL(2,q)$ does not occur.

If $G\cong \PSL(2,q)$ with $q$ odd, then
$$
|G|=\ha(q^3-q),\ |G_{P}|=\ha(q^2-q),\ |G_{P}^{(1)}|=q.
$$
The previous argument depending on the higher ramification groups
at $P$ still works as $G_P$ has two conjugacy classes of elements
of order $p$, and each of them generates $G_P^{(1)}$. Therefore,
$k=2$ and hence $g=\ha\,(q-1)$. But to show that this case cannot
actually occur, more is needed.

The $\K$--automorphism group of a hyperelliptic curve in odd
characteristic contains a central involution, say $\alpha$. Since
$\alpha$ commutes with every $p$-element of $G$, from (*) it
follows that $\alpha$ must fix $\Omega$ pointwise. If $P\in
\Omega$, then $|G_P|$ is even. But then the automorphism group
generated by $G_P$ together with $\alpha$ is not cyclic although
it fixes $P$; a contradiction.

Therefore, $\cX$ is not hyperelliptic. So, $\cX$ may be assumed to
be the canonical curve of $\K(\cX)$ embedded in $\PG(g-1,\K)$.
Then $G$ is isomorphic to a linear collineation group $\Gamma$ of
$\PG(g-1,\K)$ preserving $\cX$ such that the restriction of the
action of $\Gamma$ on $\cX$ is $G$. To simplify notation, the
symbol $G$ is used to indicate $\Gamma$, too.

\begin{lemma}
\label{psl1} Let $G\cong \PSL(2,q)$ with $q$ odd. Then
    \begin{enumerate}
    \item[(1)] $j_{g-1}(P)=2g-2$;
    \item[(2)] Let $H_{g-2}$ be a hyperplane of $\PG(g-1,\K)$ such
    that $$I(P,\cX\cap H_{g-2})=j_{g-2}(P).$$ If $H_{g-2}$ contains a
    point $R\in \Omega$ distinct from $P$, then $$I(R,\cX\cap
    H_{g-2})=2g-2-j_{g-2}(P)$$ and hence $P$ and $R$ are the only common points of $\cX$ and
    $H_{g-2}$.
    \end{enumerate}
\end{lemma}
\begin{proof} To show (1) assume on the contrary that the
osculating hyperplane $H_{g-1}$ to $\cX$ at $P$ contains a point
$S\in \cX$ distinct from $P$. Since $G_P^{(1)}$ preserves
$H_{g-1}$, the $G_P^{(1)}$-orbit of $S$ lies in $H_{g-1}$. Since
such a $G_P^{(1)}$-orbit is long, this implies that $H_{g-1}$
contains from $\cX$ at least $q$ points other than $P$. Hence,
$\deg \cX\geq j_{g-1}(P)+q$. On the other hand, $j_{g-1}(P)\geq
g-1$. Therefore,
$$\deg \cX\geq g-1+q>g-1+pg/(p-1)>2g-2,$$ contradicting\, $\deg
\cX=2g-2$.

Similar argument may be used to show (2). Again, assume on the
contrary that $H_{g-2}$ contains a point $T\in \cX$ other than $P$
and $R$. As $H_{g-1}$ does not contain $R$, $H_{g-2}$ is the
unique hyperplane through $R$ whose intersection number with $\cX$
at $P$ is $j_{g-2}(P)$. In particular, the stabiliser $H$ of $R$
in $G_P$ preserves $H_{g-2}$. Since the $H$-orbit of $T$ is long,
$H_{g-2}$ contains from $\cX$ at least $\ha\,(q-1)$ points other
than $P$ and $R$. On other hand $j_{g-2}(P)\geq g-2$. Therefore,
$$2g-2=\deg \cX\geq g-2+1+\ha\,(q-1),$$
whence $g\geq \ha\,(q+1)$, a contradiction.
\end{proof}
Since $\Omega$ is a $G_P^{(1)}$-orbit, Lemma \ref{psl1}(2) shows
that for every $R\in \Omega\setminus \{P\}$ there exists a
hyperplane $H_{g-2}(R)$ such that
$$I(P,\cX\cap H_{g-2}(R))=j_{g-2}(P),\qquad I(R,\cX\cap H_{g-2}(R))=n=2g-2-j_{g-2}(P).$$
Since such hyperplanes are distinct for distinct points $R$, from
(\ref{gnr}) it follows that
$2g-2\geq -(n+1)+q(n-1)$ with $n\geq 3$.
As $g=\ha\,(q-1)$, this leaves just one possibility, namely
$q=5,\,g=2,\,n=3$. But then $\cX$ would be hyperelliptic, a
contradiction.
Therefore the case $G\cong \PSL(2,q)$ does not occur.

If $G\cong \PSU(3,n)$ with $q=n^3$ then
$$
|G|=(n^3+1)n^3(n^2-1)/\mu,\ |G_{P}|=n^3(n^2-1)/\mu,\
|G_{P}^{(1)}|=n^3,
$$
where $\mu={\rm{gcd}}(3,n+1)$. By (\ref{Sticht13.11tris}),
$$
2g=\frac{(n^3+1)(n^2-1)} {\mu |G_{Q}|}-(n^3+1).
$$
Since (\ref{eq0}) is assumed, Lemma \ref{backgr} (i) together with
the first two assertions in Lemma \ref{GQfpf} ensure the existence
of a divisor $t\geq 1$ of $(n^2-n+1)/\mu$ such that
$|G_{Q}|=(n^2-n+1)/(t\mu)$. Hence
\begin{equation}
\label{psubis} 2g=(n-1)(t(n+1)^2-(n^2+n+1)).
\end{equation}
Since $t$ is odd, this and (\ref{eq0}) imply that $t=1$. Then
$g=\ha n(n-1)$. Therefore, (III) holds. Since $\PSU(3,n)$ is a
subgroup of $\PGU(3,n)$ of index $\mu$, the above argument works
for $G\cong \PGU(3,n)$.

If $G\cong \Sz(n)$ with $q=n^2$ and
$n=2n_0^2$ for a power $n_0\geq 2$ of $2$, then
$$
|G|=(n^2+1)n^2(n-1),\ |G_{P}|=n^2(n-1),\ |G_{P}^{(1)}|=n^2.
$$
By (\ref{Sticht13.11tris}),
$$
2g=\frac{(n+2n_0+1)(n-2n_0+1)(n-1)} {|G_{Q}|}-(n^2-1).
$$
{}From  the preceding argument depending on (\ref{eq0}) and Lemmas
\ref{backgr} and \ref{GQfpf}, there is an odd integer $t$ such
that either (A) or (B) holds, where
\begin{eqnarray*}
\mbox{(A)}\hspace*{8mm} & & 2g\,=\,(t-1)(n^2-1)-2tn_0(n-1), \quad
|G_{Q}|=(n +2n_0+1)/t;\\
\mbox{(B)}\hspace*{8mm} & & 2g\,=\, (t-1)(n^2-1)+2tn_0(n-1),\quad
|G_{Q}|=(n-2n_0+1)/t.
\end{eqnarray*}
In case (A), $t$ must be at least $3$. But then (\ref{eq0}) does
not hold
except for $n_0=2$. In the latter case, however, $t$ does not
divide $n+n_0+1=11$.
In case (B), (\ref{eq0}) implies that $t=1$. Then $g=n_0(n-1)$.
Therefore, (IV) holds.

If $G\cong \Ree(n)$ with $q=n^3$ and $n=3n_0^2$ for a power
$n_0\geq 0$ of $3$, then
$$
|G|=(n^3+1)n^3(n-1),\quad |G_{P}|=n^3(n-1), \quad
|G_{P}^{(1)}|=n^3.
$$
By (\ref{Sticht13.11tris}),
$$
2g=(n-1)\left(\frac{(n+3n_0+1)(n-3n_0+1)(n+1)}{|G_{Q}|}-(n^2+n+1)\right).
$$
Again the previous argument based on (\ref{eq0}) and Lemmas
\ref{backgr} and \ref{GQfpf} works showing this time the existence
of an integer $t\geq 1$ such that  either (A) or (B), or (C)
holds, where
\begin{eqnarray*}
\mbox{(A)} & 2g\,=\,(n-1)[(t-1)(n^2+1)-(t+1)n],\,\,
|G_{Q}|=(n +1)/t;\\
\mbox{(B)} & 2g\,=\,(n-1)[t(n^2-1)-3tnn_0+n(2t-1)-3tn_0+t-1], \\
& |G_{Q}|=(n +3n_0+1)/t;\\
\mbox{(C)} & 2g\,=\, n-1)[t(n^2-1)+3tnn_0+n(2t-1)+3tn_0+t-1], \\
& \quad |G_{Q}|=(n-3n_0+1)/t.
\end{eqnarray*}

In case (A), hypothesis (\ref{eq0}) yields that $t=2$. Then
$|G_Q|=\ha(n+1)$, hence $|G_Q|$ is even. But this is impossible as
every involution in $\Ree(n)$ has a fixed point.

If case (B) occurs, then $t\geq 2$. Since $t$ is odd and $t\neq
3$, $t$ must be at least $5$. But then (\ref{eq0}) does not hold.

If case (C) holds with $t=1$, this (C) reads
$2g=3n_0(n-1)(n+n_0+1),$ and (V) for $n>3$ follows. Otherwise,
$t\geq 5$ contradicting (\ref{eq0}).

It remains to investigate the possibility of the permutation
representation $\bar{G}$ of $G$ on $\gO$ having non-trivial
kernel. Such a kernel $M$ is a cyclic normal subgroup of $G$ whose
order is relatively prime to $p$. By Lemma \ref{GQfpf} no point
outside $\gO$ is fixed by a non-trivial element in $M$. Let
$\tilde{g}$ be the genus of the quotient curve $\cY=\cX/M$. {}From
(\ref{eq2}) applied to $M$,
$$
2g-2=|M|(2\tilde{g}-2)+(|M|-1)(q+1).
$$
By (\ref{eq0}), this implies that either $|M|=2$, or $|M|=3$ and
$\tilde{g}=0$.

Suppose first that $|M|=2,\,\tilde{g}=0$. Then, $g=\ha\,(q-1)$ and
$\cY$ is rational. So, $\cY$ may be assumed to be the projective
line $\ell$ over $\K$. Let $\Omega'$ be the set of all points of
$\ell$ which lie under the points of $\Omega$. Then
$|\Omega'|=|\Omega|$ and $\bar{G}$ acts on $\Omega'$ and $\Omega$
in the same way. So, $\bar{G}$ may be viewed as a subgroup of
$\PGL(1,\K)$ acting on a subset $\Omega'$ of $\ell$. This shows
that no non-trivial element of $\bar{G}$ fixes three distinct
points.


Assume that $\bar{G}$ has a regular normal subgroup. Arguing as in
Case (i), this yields $q+1=2^k$ with both $q$ and $k$ primes. In
particular, $q=p$ with $p-1>g=\ha\,(p-1)$. From Roquette's theorem
\cite{roquette1970}, $|G|<84(g-1)=42(p-3)$. On the other hand
$$|G|\geq 2(p+1)p,$$ as $\bar{G}$ is doubly transitive on
$\Omega$. This together with $p+1=2^k$ leaves only one
possibility, namely $p=7,\,k=3,\, g=3,\,|G|=112$ and $\bar{G}$ is
sharply $2$-transitive on $\Omega$. In particular,
$|G_P|=14,\,|G_P^{(1)}|=7$. But then (\ref{Sticht13.11tris})
yields that $|G_Q|=1$, a contradiction.

If $\bar{G}$ has no regular normal subgroup, from the
classification of Zassenhaus groups either $\bar{G}$ is a sharply
$3$-transitive group on $\Omega$, or $\bar{G}$ is $\PSL(2,q)$,
Therefore, $|G|=2q(q-1)(q+1)$ in the former case, and
$|G|=q(q-1)(q+1)$ in the latter case. In both cases,
$|G|>8(\ha\,(q-1)^3)=8g^3.$ By Theorem \ref{hennmainbound}, (II)
holds.

Suppose next that $|M|=2,\,\tilde{g}=1$. Since $G_P^{(1)}$ may be
viewed as a sub- group of $\aut(\cY)$ of the elliptic curve
$\cY=\cX/M$, the order of $G_P^{(1)}$ does not exceed $24$. Thus,
$q$ is one of the integers $2,3,4,8$. This leaves just one case,
namely $q=p=g=3$, but then (\ref{eq0}) fails.

Suppose now that $|M|=2,\,\tilde{g}\geq 2$. Then (\ref{eq0}) holds
for $\cY=\cX/M$ with $\bar{G}$ acting faithfully on
$\bar{\Omega}$. From what proven before, either $\bar{G}$ contains
a subgroup $\bar{G}'\cong \PGU(3,n)$ of index $\gcd(3,n+1)$, or
$\bar{G}\cong \Ree(n)$ with $n>3$.

In the former case, let $G'$ be the subgroup of $G$ for which
$G'/M=\bar{G'}$. Since $|M|=2$ and hence $M\subset Z(G)$, Lemma
\ref{schur} implies that $G'=M\times U$ with $U\cong \PSU(3,n)$.
Also, $U$ acts on $\Omega$ as $\PSU(3,n)$ in its natural
$2$-transitive permutation representation. Since the one-point
stabiliser $U_P$ of $U$ with $P\in \Omega$ contains a cyclic
subgroup $V$ of even order, it turns out that the the subgroup of
$G$ generated $V$ and $M$ is a not cyclic, although it fixes $P$;
a contradiction.

In the latter case, the same argument works for $G=G'=\Ree(n)$.

Finally, suppose that $|M|=3,\,\tilde{g}=0$. Then $g=q-1$. This
together with (\ref{eq0}) imply that $p>g+1=q$, a contradiction.

    \end{document}